\documentclass[12pt]{article}
\begin{document}
\title{Sieves for Twin Primes in Class I}
\author{H. J. Weber\\Department of Physics\\
University of Virginia\\Charlottesville, VA 22904\\USA}
\maketitle
\begin{abstract}
Sieves are constructed for twin primes in class I, 
which are of the form $2m\pm D,~D\geq 3$ odd. They 
are characterized by their twin-D-I rank $m.$ They 
have no parity problem. Non-rank numbers are identified 
and counted using odd primes $p\geq 5.$ Twin-D-I ranks 
and non-ranks make up the set of positive integers. 
Regularities of non-ranks allow obtaining the 
number of twin-D-I ranks. It involves considerable 
cancellations so that the asymptotic form of its 
main term collapses to the expected form, but its 
coefficient depends on $D.$           
\end{abstract}
\vspace{3ex}
\leftline{MSC: 11A41, 11N05}
\leftline{Keywords: Twin-D-I rank, non-ranks, sieve} 


\section{Introduction}

Sieve theory has developed over almost a century 
into a versatile tool of number theory~\cite{hr},
\cite{rm},\cite{hri},\cite{fi}. For twin primes 
it is the method of choice. The first genuine 
pair sieve constructed in Ref.~\cite{adhjw} for 
ordinary twin primes is adapted to twin primes 
at distance $2D\geq 6$ with $D$ odd and fixed 
throughout, except for examples. Their arithmetic 
is fairly different from distance $2$ (or 
$4$)~\cite{adhjw},\cite{hjw4}, because the 
half-distance $D$ has at least one odd prime 
divisor, whereas for ordinary prime twins of 
the form $6m\pm 1$ it has none.    

Prime numbers $p\geq 5$ are well known to be of the 
form~\cite{hw} $6m\pm 1.$ Since $2, 3$ are not of 
the form $6m\pm 1,$ they are excluded as primes in 
the following. An ordinary twin prime occurs when 
both $6m\pm 1$ are prime. Twin primes at distance 
$2D$ can be written similarly as $2m\pm D,~D$ odd 
being in class I of the classification~\cite{hjw},
\cite{hjw1} of all twin primes, the same class as 
ordinary twins of the form $2(3m)\pm 1$. 

{\bf Definition~1.1.} The base set of the sieve 
consists of all positive integers; it is \
partitioned into twin-D-I ranks and non-ranks. 
A number $m$ is called {\it twin-D-I rank} if 
$2m\pm D$ are both prime. If $2m\pm D$ are not 
both prime, then $m$ is a {\it non-rank}. 
Multiples $nq$ of divisors $q\mid D$ are trivial 
non-ranks because $2nq\pm D$ are never prime.  
 
{\bf Example~1.2.} Twin-D-I ranks for $D=3$ are 
$4, 5, 7, 8, 10, \ldots ;$ for $D=5$ they are 
$3, 4, 6, 9, 11, 12\ldots .$ Non-ranks for $D=3$ 
are $6, 9, 11, 12, 13,\ldots ;$ for $D=5$ they 
are $5, 7, 8, 10, 13, 14,\ldots .$  

Only non-ranks have sufficient regularity and 
abundance allowing us to determine the number 
of twin-D-I ranks. Therefore, our main focus 
is on non-ranks, their symmetries and abundance.      

In Sect.~2 the twin-D-I prime sieve is constructed 
based on non-ranks. In Sect.~3 non-ranks are 
identified in terms of their main properties and 
then, in Sect.~4, they are counted. In Sect.~5 
twin-D-I ranks are isolated and counted. 
Conclusions are summarized and discussed in Sect.~6.   

\section{Twin Ranks, Non-Ranks and Sieve}

It is our goal here to construct twin-D-I prime 
sieves in detail. We need the following arithmetical 
function~\cite{adhjw},\cite{hjw4}. 

{\bf Definition~2.1.} Let $x$ be real. Then $N(x)$ is the 
integer nearest to $x.$ The ambiguity for $x=n+\frac{1}{2}$ 
with integral $n$ will not arise in the following.    

{\bf Lemma~2.2.} {\it Let $p\geq 5$ be prime. Then} 
\begin{eqnarray}
N(\frac{p}{6})=\{\begin{array}{ll}\frac{p-1}{6},~\rm{if}~
p\equiv 1\pmod{6};\\\frac{p+1}{6},~\rm{if}~p\equiv 
-1\pmod{6}.\\\end{array}
\end{eqnarray}
{\it Proof.} This is obvious from Def.~2.1 by 
substituting $p=6m\pm 1.~\diamond$ 

{\bf Lemma~2.3} {\it Let $p\geq 5$ be prime and 
$(p,D)=1$. Then the numbers   
\begin{eqnarray}\nonumber 
&&k(n,p)^+=np+3DN(\frac{p}{6}),~n=0,1,2,\ldots\\
&&k(n,p)^-=np-3DN(\frac{p}{6})>0,~n>\frac{D+1}{2} 
\end{eqnarray}
are non-ranks. There are $2=2^{\nu(p)}$ (single) 
non-rank progressions to the prime $p.$  
  
(a) If $p\equiv 1\pmod{6}$ the non-rank 
$k(n,p)^+$ generates the pair}  
\begin{eqnarray}
2k(n,p)^+\pm D=((2n+D)p-2D,(2n+D)p),
\label{p++}
\end{eqnarray} 
{\it and the non-rank $k(n,p)^-$ the pair}  
\begin{eqnarray}
2k(n,p)^-\pm D=((2n-D)p,(2n-D)p+2D),~2n>D+1.
\label{p+-}
\end{eqnarray} 
{\it (b) If $p\equiv -1\pmod{6}$ the non-rank 
$k(n,p)^+$ generates the pair}  
\begin{eqnarray}
2k(n,p)^+\pm D=((2n+D)p,(2n+D)p+2D);
\label{p-+}
\end{eqnarray} 
{\it and the non-rank $k(n,p)^-$ the pair} 
\begin{eqnarray}
2k(n,p)^-\pm D=((2n-D)p-2D,(2n-D)p),~2n>D+1.
\label{p--}
\end{eqnarray}
{\it All pairs contain a composite number.}

Clearly, all these non-ranks are symmetrically distributed 
at equal distances $3D N(p/6)$ from multiples of each prime 
$p\geq 5,$ except for prime divisors of $D.$  
  
{\it Proof.} Let $p\equiv 1\pmod{6}$ be prime and $n\geq 0$ 
an integer. Then $2k(n,p)^+\pm D=2np+D(p-1)\pm D$ 
by Lemma~2.2 and $2k^+$ is sandwiched by the pair in 
Eq.~(\ref{p++}) which contains a composite number. Hence 
$k(n,p)^+$ is a non-rank. For $2n>D+1,$ the same happens 
in Eq.~(\ref{p+-}), so $k^-$ is a non-rank. 

If $p\equiv -1\pmod{6}$ and prime, then $2k(n,p)^+\pm 
D=2np+D(p+1)\pm D$ by Lemma~2.2 and $k^+$ leads to the 
pair in Eq.~(\ref{p-+}) which contains a composite 
number again. For $2n>D+1,$ the same happens in 
Eq.~(\ref{p--}), so $k^-$ is a non-rank.~$\diamond$

The $k(n,p)^{\pm}$ yield pairs $2k^{\pm}\pm D$ with one 
or two composite entries that are twin-D-I prime analogs 
of multiples $np,~n>1,$ of a prime $p$ in Eratosthenes' 
prime sieve~\cite{hw}. Non-ranks form the sieving set. 

The converse of Lemma~2.3 holds, i.e. nontrivial non-ranks 
are organized in terms of arithmetic progressions with 
primes $\geq 5$ (and their products) as periods. This 
makes it a cornerstone of the pair sieves.    
 
{\bf Lemma~2.4.} {\it If $k$ is a nontrivial non-rank, 
there is a prime $p\geq 5$ and an integer $\lambda$ 
so that $k=k(\lambda,p)^+$ or} $k=k(\lambda,p)^-.$ 

{\it Proof.} If $k\equiv 0\pmod{3}$ and $D\equiv 
1\pmod{6},$ then $2k+D\equiv 1\pmod{6}.$ Let 
$2k+D=pK$ be composite, where $p\geq 5$ is the 
smallest prime divisor $\not| D$. Then $2k+D\neq 
3^\nu,~\nu\geq 1$ obviously. If $p\equiv 1\pmod{6}$ 
then $K\equiv 1\pmod{6}.$ So $p=6m+1,~K=6\kappa+1, 
k=3k',~D=6d+1$ and  
\begin{eqnarray}
2k+D=6^2 m\kappa+6(m+\kappa)+1,~k'+d=6m\kappa
+m+\kappa=p\kappa+\frac{p-1}{6}. 
\end{eqnarray}
Hence 
\begin{eqnarray}
k=3p\kappa-3d+3\frac{p-1}{6}=p(3\kappa-
\frac{D-1}{2})+3DN(\frac{p}{6})
\end{eqnarray} 
and $\lambda=3\kappa-\frac{D-1}{2}.$

If $p\equiv -1\pmod{6}$ then $K\equiv -1\pmod{6},$ 
i.e. $p=6m-1,~K=6\kappa-1$ and   
\begin{eqnarray}
2k+D=6^2 m\kappa-6(m+\kappa)+1,~k'+d=6m\kappa-(m
+\kappa)=p\kappa-\frac{p+1}{6},
\end{eqnarray}
then 
\begin{eqnarray}
k=3p\kappa-3d-3\frac{p+1}{6}=p(3\kappa+
\frac{D-1}{2})-3DN(\frac{p}{6}),
\end{eqnarray} 
and $\lambda=3\kappa+\frac{D-1}{2}.$

If $2k-D=pK$ i.e. is composite and $p=6m+1,$ 
then $K=6\kappa-1$ because $2k-D=6k'-6d-1.$ 
Hence  
\begin{eqnarray}
2k-D=6^2 m\kappa+6(\kappa-m)-1,~k'-d=p\kappa
-\frac{p-1}{6} 
\end{eqnarray}
and
\begin{eqnarray}
k=p(3\kappa+\frac{D-1}{2})-3DN(\frac{p}{6}). 
\end{eqnarray} 
So $\lambda=3\kappa+\frac{D-1}{2}.$

If $p=6m-1$ then $K=6\kappa+1$ and $2k-D=6(k'-\kappa)-1,$ 
\begin{eqnarray}
2k-D=6^2 m\kappa+6(m+\kappa)-1,~k'+d=6m\kappa-(m+\kappa)
=p\kappa-\frac{p+1}{6}.
\end{eqnarray}
Hence
\begin{eqnarray}
k=3p\kappa+3d-3\frac{p+1}{6}=p(3\kappa
-\frac{D-1}{2})+3DN(\frac{p}{6})
\end{eqnarray} 
and $\lambda=3\kappa-\frac{D-1}{2}.$

If $D\equiv 3\pmod{6},$ then $2k\pm D=3(2k'+2d+1)$ 
is always composite. These trivial non-rank cases 
are ignored in the following, except when counting 
non-ranks.  

If $D=6d-1\equiv -1\pmod{6},$ then there are 
four cases for $2k+D$ or $2k-D$ composite, combined 
with the options for $p$ and $K$ as in (i), which 
are all handled the same way and then lead to 
similar results.  

(ii) If $k=3k'+1\equiv 1\pmod{3}$ and $D=6d-1\equiv 
-1\pmod{6},$ then $2k+D=6k'+6d+1\equiv 1\pmod{6}.$  
Let $2k+D=pK$ be composite and $p=6m+1.$ Then 
$K=6\kappa+1$ and 
\begin{eqnarray}
2k+D=p\cdot K=(6m+1)(6\kappa+1)=6^2m\kappa
+6(m+\kappa)+1.
\end{eqnarray}       
Hence 
\begin{eqnarray}\nonumber
k'+d&=&6m\kappa+m+\kappa=p\kappa+\frac{p-1}{6},\\
k&=&3k'+1=[3\kappa-\frac{D+1}{2}]p+3D N(\frac{p}
{6}),
\end{eqnarray}
so $\lambda=3\kappa-\frac{D+1}{2}.$ 

If $D=6d+3,$ then $2k+D=6k'+5+6d\equiv -1\pmod{6}
=pK.$ If $p=6m+1,$ then $K=6\kappa-1$ and  
\begin{eqnarray}\nonumber
2k+D&=&6^2 m\kappa+6(\kappa-m)-1,\\\nonumber
k'+d+1&=&p\kappa-\frac{p-1}{6},\\\nonumber
k&=&3p\kappa-3(d+1)+1-\frac{p-1}{2}\\
&=&p(3\kappa-\frac{D-3}{2}-2)+3DN(\frac{p}{6}),
\end{eqnarray}
so $\lambda=3\kappa-\frac{D-3}{2}-2$. The case 
where $p=6m-1,~K=6\kappa+1$ is handled similarly. 
All other cases lead to trivial non-ranks.  

(iii) The cases for $k=3k'-1$ are similar and 
handled in the same way.$~\diamond$ 

{\bf Theorem~2.5. (Prime Pair Sieves)} {\it Let} ${\cal 
P}=\{(2m-D\geq 3,2m+D): 2m\not\equiv 0\pmod{q}, q\mid 
D,~m~\rm{integral}\}$ {\it be the set of pairs with 
entries $\geq 3$ of natural numbers at distance $2D,~D$ 
odd. Upon striking all pairs identified by non-ranks of 
Lemma~2.3, only (and all) twin-D-I prime pairs are left.} 

Since after sieving only twin-D-I ranks are left that 
lead to prime pairs at distance $2D$ (and no composites) 
the sieves have no parity problem.     
 
{\it Proof.} For $2m-D\geq 3$ divide $2m\pm D$ by all 
primes $p<\sqrt{2m+D}.$ Then $m$ is a non-rank if there 
is a prime $p$ such that $(2m-D)/p$ or $(2m+D)/p$ (or 
both) is integral. All such $m$ are struck from the set 
of positive integers. Then all remaining integers are 
twin-D-I ranks.~$\diamond$    

\section{Identifying Non-Ranks} 

Here it is our goal to characterize and systematically  
identify non-ranks among natural numbers. 

{\bf Definition~3.1} Let $p\geq 5$ be the minimal 
prime of a non-rank. Then $p$ is its parent prime.  

The non-ranks to parent prime $5$ are, by Lemma~2.3,   
\begin{eqnarray}
k^+=5n+3D,n>-[\frac{3D}{5}];k^-=5n-3D,5n>3D;(n,D)=1,
5\not|D,  
\end{eqnarray} 
where $[x]$ is the largest integer below $x,$ as 
usual. These $k^{\pm}$ form the set ${\cal A}_5^-
={\cal A}_5.$ 
  
Note that $5$ is the most effective non-rank 
generating prime number (except when $5|D$). If it 
were excluded like $3$ then many numbers would be 
missed as non-ranks. 

In contrast to ordinary twin primes~\cite{adhjw} 
the arithmetic function values $N(p'/6), N(p/6)$ 
do not suffice to characterize twin-D-I primes 
$p'=p+2D.$ 

{\bf Lemma~3.2.} {\it Let $p'>p$ be primes.  
Then $p'=p+2$ are ordinary prime twins iff} 
$N(\frac{p'}{6})=N(\frac{p}{6}).$  

{\it Proof.} See Theor.~3.6 of Ref.~\cite{adhjw}. 

Lemma~3.2 generalizes to $D\geq 3$ as follows. 

{\bf Corollary~3.3.} {\it Let $p'>p\geq 5$ be 
primes with $p'\equiv p\pmod{6}.$  

(i) If $D\equiv 0\pmod{3}$ then $p'=p+2D$ holds 
iff $N(\frac{p'}{6})=\frac{D}{3}+N(\frac{p}{6}).$

(ii) If $D\equiv 1\pmod{3}$ then $p'=p+2(D-1)$  
iff $N(\frac{p'}{6})=\frac{D-1}{3}+N(\frac{p}{6}).$

(iii) If $D\equiv -1\pmod{3}$ then $p'=p+2(D+1)$  
iff} $N(\frac{p'}{6})=\frac{D+1}{3}+N(\frac{p}{6}).$

{\it Proof.} (i) $\frac{p'\mp 1}{6}=\frac{D}{3}+
\frac{p\mp 1}{6}$ is equivalent to $p'=p+2D.$ 
(ii) $\frac{p'\mp 1}{6}=\frac{D-1}{3}+
\frac{p\mp 1}{6}$ is equivalent to $p'=p+2(D-1).$
(iii) $\frac{p'\mp 1}{6}=\frac{D+1}{3}+
\frac{p\mp 1}{6}$ is equivalent to $p'=p+2
(D+1).~\diamond$ 

It is straightforward to relax the constraint 
$p'\equiv p\pmod{6}$ to include $p'\equiv p\pm 
2\pmod{6}.$ 

We now consider systematically common (or double) 
non-ranks of pairs of primes. We start with 
ordinary twin primes.   

{\bf Theorem~3.4.} {\it Let $p'>p\geq 5$ and 
$N(\frac{p'}{6})=N(\frac{p}{6}).$ Then (i)  
\begin{eqnarray}
pp'n\pm 3DN(\frac{p}{6})=p'pn\pm 3DN(\frac{p'}{6}) 
\label{otp}
\end{eqnarray}
are two common non-rank progressions of $p$ and $p'.$
(ii) If $r, r'$ solve 
\begin{eqnarray}\nonumber 
&&(r'-r\pm D)p=2r\pm D,~p\equiv 1\pmod{6}\\
&&(r'-r\pm D)p=2r\mp D,~p\equiv -1\pmod{6}, 
\label{otpo}
\end{eqnarray}
then 
\begin{eqnarray}
p(p'n+r')\pm 3DN(\frac{p}{6})=p'(pn+r)\mp 3DN(\frac{p'}{6}) 
\label{otpo1}
\end{eqnarray}
are the other two common non-rank progressions of $p$ and} 
$p'.$

{\it Proof.} By Lemma~3.2, $p'=p+2$ and $p$ are 
ordinary twin primes and Eq.~(\ref{otp}) is valid 
obviously, with the lhs a non-rank to $p$ and the 
rhs a non-rank to $p'.$ (ii) If $2r=\mp D+\lambda p, 
r'=r+\lambda\mp D$ for odd $\lambda$ so that $-p<2r<p, 
-p'<2r'<p',$ solving Eq.~(\ref{otpo}) for $p\equiv 
1\pmod{6},$ then Eq.~(\ref{otpo1}) is verified to be 
equivalent to Eq.~(\ref{otpo}), its lhs being a 
non-rank to $p$ and rhs a non-rank to $p'.$ For $
p\equiv -1\pmod{6}$ in Eqs.~(\ref{otpo}),(\ref{otpo1}) 
the cases are treated similarly.~$\diamond$ 
   
{\bf Theorem~3.5.} {\it Let $p'>p\geq 5$ 
be primes with $(p,D)=1=(p',D)$. (i) If 
$p'\equiv p\pmod{6},$ then $p'=p+6l,~l\geq 
1,~ N(\frac{p'}{6})=N(\frac{p}{6})+l$ and two 
common non-rank progressions of $p',p$ are} 
\begin{eqnarray}
p[p'n+r']\pm 3DN(\frac{p}{6})=
p'[pn+r]\pm 3DN(\frac{p'}{6})
\label{snr1}
\end{eqnarray}
{\it provided $r, r'$ solve} 
\begin{eqnarray}
(r'-r)p=3l(2r\pm D). 
\label{snr11}
\end{eqnarray}
{\it The solution of Eq.~(\ref{snr11}), $2r=\mp 
D+p\lambda_{\pm},~-p<2r<p$ for odd $\lambda_{\pm}$ 
with $r'=r+\lambda_{\pm},~-p'<2r'<p'$ on the 
lhs of Eq.~(\ref{snr1}) yields a non-rank to $p$ 
and, on the rhs, a non-rank to} $p'.$     

{\it If $r, r'$ solve} 
\begin{eqnarray}\nonumber
&&(r'-r\pm D)p=6lr\mp D(3l-1),~p\equiv 1\pmod{6}\\
&&(r'-r\pm D)p=6lr\mp D(3l+1),~p\equiv -1\pmod{6},    
\label{snr12}
\end{eqnarray}
{\it then two more common non-rank progressions are}
\begin{eqnarray}
p[p'n+r']\pm 3DN(\frac{p}{6})=
p'[pn+r]\mp 3DN(\frac{p'}{6}).
\label{snr13}
\end{eqnarray}
{\it (ii) If $p'\equiv 1\pmod{6},~p\equiv 
-1\pmod{6}$ then $p'=p+6l+2,~l\geq 0,
~N(\frac{p'}{6})=N(\frac{p}{6})+l,$ and two  
common non-rank progressions of $p',p$ are}
\begin{eqnarray}
p[p'n+r']\pm 3DN(\frac{p}{6})=
p'[pn+r]\pm 3DN(\frac{p'}{6})
\label{snr3}
\end{eqnarray}
{\it provided}
\begin{eqnarray}
(r'-r)p=2r(3l+1)\pm 3lD. 
\label{snr33}
\end{eqnarray}
{\it If $l=0$ then $r'=r=0;$ see Theor.~3.4. 
For $l\geq 1,~2r(3l+1)=\mp 3Dl+p\lambda, r'=r+\lambda$ 
solve Eq.~(\ref{snr33}). There is a unique pair $r',r$ 
with $-p<2r<p,-p'<2r'<p'.$}   

{\it If $r, r'$ solve} 
\begin{eqnarray}
(r'-r\pm D)p=3l(2r\mp D)+2r,   
\label{snr31}
\end{eqnarray}
{\it then two more common non-rank progressions are}
\begin{eqnarray}
p[p'n+r']\pm 3DN(\frac{p}{6})=
p'[pn+r]\mp 3DN(\frac{p'}{6}).
\label{snr32}
\end{eqnarray}
{\it For appropriate $\lambda,$ the solution $r'=r+
\lambda\mp D, 3l(2r\mp D)+2r=p\lambda$ is unique.}

{\it (iii) If $p'\equiv -1\pmod{6},~p\equiv 
1\pmod{6}$ then $p'=p+6l-2,~l\geq 1,
~N(\frac{p'}{6})=N(\frac{p}{6})+l,$ and two 
common non-rank progressions of $p',p$ are}
\begin{eqnarray}
p[p'n+r']\pm 3DN(\frac{p}{6})=
p'[pn+r]\pm 3DN(\frac{p'}{6})
\label{snr4}
\end{eqnarray}
{\it provided}
\begin{eqnarray}
(r'-r)p=2r(3l-1)\pm 3Dl. 
\label{snr44}
\end{eqnarray}
{\it Again, for appropriate $\lambda$ the solution 
$r'=r+\lambda,2r(3l-1)=p\lambda\mp 3Dl$ is unique.}

{\it If $r, r'$ solve} 
\begin{eqnarray}
(r'-r\pm D)p=2r(3l-1)\mp 3lD,  
\label{snr41}
\end{eqnarray}
{\it then two more common non-rank progressions are}
\begin{eqnarray}
p[p'n+r']\pm 3DN(\frac{p}{6})=
p'[pn+r]\mp 3DN(\frac{p'}{6}).
\label{snr42}
\end{eqnarray}
{\it The solution $2r(3l-1)=\pm 3lD+p\lambda, r'=r+
\lambda\mp D$ is unique for appropriate} $\lambda.$

Note that there are $4=2^{\nu(pp')}$ 
arithmetic progressions of common or double 
non-ranks to the primes $p', p$ in all cases. 

{\it Proof.} By substituting $p',N(p'/6)$ 
in terms of $p,N(p/6)$ and $l,$ respectively, 
it is readily verified that Eqs.~(\ref{snr1}), 
(\ref{snr11}) are equivalent, as are 
(\ref{snr12}), (\ref{snr13}), and (\ref{snr31}), 
(\ref{snr32}), and (\ref{snr4}), (\ref{snr44}), 
and (\ref{snr41}), (\ref{snr42}). As in (i) 
there is a unique solution $(r, r')$ in all 
other cases as well.~$\diamond$ 

{\bf Theorem~3.6. (Triple non-ranks)} {\it Let 
$5\leq p<p'<p''$ (or $5\leq p<p''<p'$, or 
$5\leq p''<p<p'$) be different odd primes such 
that $(p,D)=1=(p',D)=(p'',D)$. Then each case 
in Theor.~3.5 of four double non-ranks leads to 
$8=2^{\nu(pp'p'')}$ triple non-ranks of $p, p', 
p''.$ At two non-ranks per prime, there are at 
most $2^3$ triple non-ranks.}  

{\it Proof.} It is based on Theor.~3.5 and 
similar for all its cases. Let's take (i) 
and substitute $n\to p''n+\nu,~0\leq 
\nu<p''$ in Eq.~(\ref{snr1}) which, upon 
dropping the term $p''p'pn,$ yields on the 
lhs 
\begin{eqnarray}
pp'\nu+pr'-3DN(\frac{p}{6})=p''\mu\pm 
3DN(\frac{p''}{6}).
\label{31}
\end{eqnarray} 
Since $(pp',p'')=1$ there is a unique residue $\nu$ 
modulo $p''$ so that the lhs of Eq.~(\ref{31}) is 
$\equiv\pm 3DN(\frac{p''}{6})\pmod{p''},$ and 
this determines $\mu.$ As each sign case leads to 
such a triple non-rank solution, it is clear that 
there are $2^3$ non-ranks to $p, p', p''.~\diamond$ 
 
{\bf Theorem~3.7. (Multiple non-ranks)} {\it Let 
$5\leq p_1<\cdots <p_m$ be $m$ different primes 
with $(p_i,D)=1$. Then there are $2^m$ arithmetic 
progressions of $m-$fold non-ranks to the primes} 
$p_1,\ldots,p_m.$ 

{\it Proof.} This is proved by induction on $m.$   
Theors.~3.5 and 3.6 are the $m=2, 3$ cases. If 
Theor.~3.7 is true for $m$ then for any case 
$5\leq p_{m+1}<p_1<\cdots<p_m,$ or $\ldots,$ 
$5\leq p_1<\cdots<p_{m+1},$ we substitute in 
an $m-$fold non-rank equation $n\to p_{m+1}n 
+\nu$ as in the proof of Theor.~3.6, again 
dropping the $n\prod_1^{m+1} p_i$ term. Then 
we get   
\begin{eqnarray}\nonumber
&&p_1(p_2(\cdots(p_m\nu+r_m)+\cdots+r_2)+
3DN(\frac{p_1}{6})\\&&=p_{m+1}\mu\pm 
3DN(\frac{p_{m+1}}{6})
\label{m}
\end{eqnarray}  
with a unique residue $\nu\pmod{p_{m+1}}$ so 
that the lhs of Eq.~(\ref{m}) becomes 
$\equiv 3DN(\frac{p_{m+1}}{6})\pmod{p_{m+1}},$ 
which then determines $\mu.$ In case the lhs of 
Eq.~(\ref{m}) has $p_1(\ldots)-3DN(p_1/6)$ the 
argument is the same. This yields an $(m+1)-$fold 
non-rank progression since each sign in 
Eq.~(\ref{m}) gives a solution. Hence there 
are $2^{m+1}$ such non-ranks. At two non-ranks 
per prime, there are at most $2^{m+1}$ non-rank 
progressions.~$\diamond$  

\section{Counting Non-Ranks}

If we subtract for case (i) in Theor.~3.5, say, 
the four common non-rank progressions 
corresponding to the solutions $-p_i<r_i<p_i,
\ldots,$ this leaves in ${\cal A}_{p'}^-=\{p'n\pm 
3D\frac{p'+1}{6}\}$ the following progressions 
$p'pn\pm 3D\frac{p'+1}{6},\ldots, p'[np+r_1]+
3D\frac{p'+1}{6},\ldots, p'[np+r_2]-3D\frac{p'+1}
{6},\ldots, p'[np+r_3]+3D\frac{p'+1}{6},\ldots, 
p'[np+r_4]-3D\frac{p'+1}{6},\ldots, 
p'np\pm 3D\frac{p'+1}{6}.$ 
 
We summarize this as follows. 

{\bf Lemma~4.1.} {\it $p'>p\geq 5$ be prime 
such that $(p,D)=1=(p',D)$. Removing the 
nontrivial common non-ranks of $p',p$ from 
the set of all non-ranks of $p'$ leaves 
arithmetic progressions of the form $p'np+l;
~n\geq 0,$ where $l>0$ are given nonnegative 
integers.}  

{\bf Proposition~4.2.} {\it Let $p\geq p'\geq 5$ 
be prime with $(p,D)=1$. Then the set of non-ranks 
to parent prime $p,~{\cal A}_p,$ is made up of 
arithmetic progressions ${\bar L}(p)n+a,~n\geq 
0$ with ${\bar L}(p)=\prod_{5\leq p'\leq p,
(p',D)=1}p$ and $a>0$ given integers.}  

{\it Proof.} Let $p=6m\pm 1.$ We start from the 
set ${\cal A}_p^\pm=\{pn\pm 3DN(\frac{p}{6})\}.$ 
Removing the non-ranks common to $p$ and $5$ by 
Theor.~3.5 leaves arithmetic progressions of the 
form $5pn+l,~n\geq 0$ where $l>0$ are given 
integers provided $5\not|D$. Continuing this 
process to the largest prime $p'<p$ leaves in 
${\cal A}_p$ arithmetic progressions of the 
form ${\bar L}(p)n+a,~n\geq 0$ with $L(p)=
\prod_{5\leq p'\leq p,(p',D)=1}p'$ and $a>0$ 
a sequence of given integers independent of 
$n.~\diamond$         

{\bf Proposition~4.3.} {\it Let $p\geq 
p'\geq 5$ be primes such that $(p,D)=1=
(p',D)$ and $G(p)$ the number of nontrivial 
non-ranks ${\bar L}(p)n+a\in {\cal A}_p$ over 
one period ${\bar L}(p)$ corresponding to 
arithmetic progressions ${\bar L}(p)n+a\in 
{\cal A}_p.$ Then} $G(p)=2\prod_{5\leq p'<p,
(p',D)=1}(p'-2).$  

Note that $G(p)<{\bar L}(p)$ both increase 
monotonically as $p\to\infty.$ 

{\it Proof.} In order to determine $G(p)$ 
we have to eliminate all non-ranks of 
primes $5\leq p'<p$ from ${\cal A}_p.$ 
As in Theor.~3.5 we start by subtracting 
the fraction $2/5$ of non-ranks to $p'=5$ 
from the interval $1\leq a\leq {\bar L}(p),$ 
then $2/7$ for $p'=7$ and so on for 
all $p'<p.$ The factor of $2$ is due to the 
symmetry of non-ranks around each multiple 
of $p'$ according to Lemma~2.3. This leaves 
$p\prod_{5\leq p'<p,(p',D)=1}(p'-2)$ numbers 
$a$. The fraction $2/p$ of these are the 
non-ranks to parent prime $p.~\diamond$   

Prop.~4.3 implies that the fraction of 
non-ranks related to a prime $p$ in the 
interval occupied by ${\cal A}_p,$  
\begin{eqnarray}
q(p)=\frac{G(p)}{{\bar L}(p)}=\frac{2}{p}
\prod_{5\leq p'<p,(p',D)=1}\frac{p'-2}{p'},
\label{qp}
\end{eqnarray}
where $p'$ is prime, decreases 
monotonically as $p$ goes up.  

{\bf Definition~4.4.} Let $p\geq p'\geq 5$ be 
prime such that $(p,D)=1=(p',D)$. The supergroup 
${\cal S}_p=\bigcup_{5\leq p'\leq p;p,(p',D)=1}
{\cal A}_{p'}$ contains the sets of non-ranks 
corresponding to arithmetic non-rank progressions 
$a+{\bar L}(p')n$ of all ${\cal A}_{p'},~5\leq 
p'\leq p;(p',D)=1=(p,D).$   

Thus, each supergroup ${\cal S}_p$ contains 
nested sets of non-ranks related to primes 
$5\leq p'\leq p,(p',D)=1=(p,D)$

Let us now count prime numbers from $p_1=5$ 
on provided $5\not|D$, omitting prime divisors 
of $D$ along with $2$ and $3.$
 
{\bf Proposition~4.5.} {\it Let $p_j\geq 5$ 
be the $j$th prime. (i) Then the number of 
nontrivial non-ranks $a\in {\cal A}_{p_i}$ 
corresponding to arithmetic progressions 
related to a prime $p_i<p_j,(p_i,D)=1=(p_j,D)$}
\begin{eqnarray}
G(p_i)=\frac{{\bar L}(p_j)}{{\bar L}(p_i)}G(p_j)
=\frac{2{\bar L}(p_j)}{p_i}\prod_{5\leq p<p_i;
(p,D)=1=(p_i,D)}\frac{p-2}{p}=q(p_i){\bar L}(p_j),
\end{eqnarray}
{\it where $p$ is prime, monotonically decreases 
as $p_i$ goes up. (ii) The number of nontrivial 
non-ranks in a supergroup ${\cal S}_{p_j}$ over 
one period ${\bar L}(p_j)$ is}
\begin{eqnarray}
S(p_j)={\bar L}(p_j)\sum_{5\leq p\leq p_j;(p,
D)=1=(p_j,D)}q(p)={\bar L}(p_j)\left(1-\prod_{
5\leq p\leq p_j;(p,D)=1=(p_j,D)}\frac{p-2}
{p}\right).
\label{lsp} 
\end{eqnarray}
{\it (iii) The fraction of non-ranks of their 
arithmetic progressions in the (first) interval 
$[1,{\bar L}(p_j)]$ occupied by the supergroup} 
${\cal S}_{p_j},$
\begin{eqnarray}
Q(p_j)=\frac{S(p_j)}{{\bar L}(p_j)}=\sum_{5\leq 
p\leq p_j;(p,D)=1=(p_j,D)}q(p)=1-\prod_{5\leq 
p\leq p_j;(p,D)=1=(p_j,D)}\frac{p-2}{p},
\label{qp1}
\end{eqnarray}
{\it increases monotonically as $p_j$ goes up.} 

{\it Proof.} (i) follows from Prop.~4.3 
and Eq.~(\ref{qp}). (ii) and (iii) are 
equivalent and are proved by induction as follows, 
using Def.~4.4 in conjunction with Eq.~(\ref{qp}).  

From Eq.~(\ref{qp}) we get $q_1=2/p_1$ which is 
the case $j=1,~p_j=5$ of Eq.~(\ref{qp1}). Assuming 
Eq.~(\ref{qp1}) for $p_j,$ we add $q_{j+1}$ of 
Eq.~(\ref{qp}) and obtain
\begin{eqnarray}\nonumber
\sum_{i=1}^{j+1}q(p_i)&=&1-\prod_{i=1}^j\frac{
p_i-2}{p_i}+\frac{2}{p_{j+1}}\prod_{i=1}^j\frac{
p_i-2}{p_i}\\&=&1-\prod_{i=1}^{j+1}\frac{p_i-2}
{p_i}.
\end{eqnarray}
The extra factor $0<(p_{j+1}-2)/p_{j+1}<1$ shows 
that $q(p_j), x(p_j)$ in Eq.~(\ref{xp}) decrease 
monotonically as $p_j\to p_{j+1}$ while $Q(p_j)$ 
increases as $j\to\infty.~\diamond$

{\bf Definition~4.6.} Since ${\bar L}(p)>S(p),$ 
there is a set ${\cal R}_p$ of {\it remnants} 
$r\in [1, {\bar L}(p)]$ such that $r\not\in{\cal 
S}_p, (r,D)=1.$   

{\bf Lemma~4.7.} {\it (i) The number $R(p_j)$ of 
remnants in a supergroup, ${\cal S}_{p_j},$ with 
$(p_j,D)=1,p_j$ prime is}   
\begin{eqnarray}\nonumber
R(p_j)&=&{\bar L}(p_j)-S(p_j)={\bar L}(p_j)
(1-Q(p_j))=\prod_{5\leq p\leq p_j;(p,D)=1
=(p_j,D)}(p-2)\\&=&\frac{1}{2}G(p_{j+1}).
\label{rps}
\end{eqnarray} 

{\it (ii) The fraction of such remnants in} 
${\cal S}_{p_j},$   
\begin{eqnarray}
x(p_j)&=&\frac{R(p_j)}{{\bar L}(p_j)}=1-Q(p_j)
=\prod_{5\leq p\leq p_j;(p,D)=1=(p_j,D)}\frac{
p-2}{p},\label{xp}
\end{eqnarray}
{\it where $p$ is prime, decreases monotonically 
as} $p_j\to\infty.$ 

{\it Proof.} (i) follows from Def.~4.6 in 
conjunction with Eq.~(\ref{lsp}) and (ii) 
from Eq.~(\ref{rps}). Equation~(\ref{rps}) 
follows from Eq.~(\ref{qp1}). 
       
\section{Remnants and Twin Ranks}

When all primes $5\leq p\leq p_j,(p_j,D)=1$ and 
appropriate nonnegative integers $n$ are used in 
Lemma~2.3 one will find all non-ranks $k<M(j+1)
\equiv (p_{j+1}^2-D^2)/2.$ By subtracting these 
non-ranks from the set of positive integers 
$N\leq M(j+1)$ all and only twin-D-I ranks 
$t<M(j+1)$ are left among the remnants provided 
trivial non-ranks are also eliminated. If a 
non-rank $k$ is left then $2k\pm D$ must have 
prime divisors that are $>p_j$ according to 
Lemma~2.3, which is impossible. 

{\bf Definition~5.1.} Let $p_j, p_{j+1}\geq 5$ 
be prime such that $(p_j,D)=1=(p_{j+1},D)$. 
Then all $t<M(j+1)=(p_{j+1}^2-D^2)/2$ in a remnant 
${\cal R}_{p_j}$ of a supergroup ${\cal S}_{p_j}$ 
are twin-D-I ranks. These twin ranks are called 
{\it front twin ranks}. 

Twin ranks are located among the remnants 
${\cal R}_p$ for any prime $p\geq 5,(p,D)=1$. 
Our goal is to determine the number of twin-D-I 
ranks.   

{\bf Theorem~5.2.} {\it Let $R_0$ be the number 
of remnants of the supergroup ${\cal S}_{p_j},$ 
where $p_j$ is the $j$th prime number with 
$(p_j,D)=1,~p_j>p,~\forall p|D$ and $M(j+1)=
(p_{j+1}^2-D^2)/2$. Then the number 
$R=\pi_2(2{\bar L}(p_j)+D)/2$ 
of twin-D-I ranks within the remnants of the 
supergroup ${\cal S}_{p_j}$ is given by} 
\begin{eqnarray}
R\prod_{p|D}(1-\frac{1}{p})^{-1}=R_0+
\sum_{(n,D)=1,p_j>n}\mu(n)2^{\nu(n)}\bigg[
\frac{{\bar L}(p_j)-M(j+1)}{n}\bigg].
\label{tr} 
\end{eqnarray}

Here ${\bar L}(p_j)=\prod_{5\leq p\leq p_j,
(p,D)=1}p,$ $R_0=\prod_{5\leq p\leq p_j,
(p,D)=1}(p-2)$ with $p$ prime, and $n$ runs 
through all products of primes $p_j<p\leq 
(2{\bar L}(p_j)+1)/D$ relatively prime to 
$D$. The upper limit $(2{\bar L}(p_j)+1)/
D$ comes about because $3DN(p/6)$ is the 
lowest possible non-rank of a prime number 
$p$ according to Lemma~2.2.     

The argument of the twin-prime counting function 
$\pi_2$ is $2{\bar L}(p_j)+D$ because, if ${\bar 
L}(p_j)$ is the last twin-D-I rank of the interval 
$[1, {\bar L}(p_j)],$ then $2{\bar L}(p_j)\pm D$ 
are the corresponding twin-D-I primes.    

{\it Proof.} According to Prop.~4.5 the 
supergroup ${\cal S}_{p_j}$ has $S(p_j)=
{\bar L}(p_j)\cdot\left(1-\prod_{5\leq p\leq p_j}
\frac{p-2}{p}\right)$ non-ranks. Subtracting 
these from the interval $[1, {\bar L}(p_j)]$ 
that the supergroup occupies gives $R_0=\prod_{
5\leq p\leq p_j,(p,D)=1}(p-2)$ for the number 
of remnants which include twin-D-I ranks and 
non-ranks to primes $p_j<p\leq (2{\bar L}(p_j)
+1)/D.$ The latter are 
\begin{eqnarray}
M(j+1)<pn\pm 3DN(\frac{p}{6})\leq {\bar L}
(p_j),~M(j+1)=(p_{j+1}^2-D^2)/2, 
\end{eqnarray} 
or 
\begin{eqnarray}
0<n\leq \frac{{\bar L}(p_j)-M(j+1)}{p}, 
\end{eqnarray}
which have to be subtracted from the remnants 
to leave just twin-D-I ranks. Correcting for 
double counting of common non-ranks to two 
primes using Theor.~3.5, of triple non-ranks 
using Theor.~3.6 and multiple non-ranks 
using Theor.~3.7 and eliminating trivial 
non-ranks to prime divisors of $D,$ which 
leads to the factor on the lhs of 
Eq.~(\ref{tr}), we obtain 
\begin{eqnarray}\nonumber
&&R\prod_{p|D}(1-\frac{1}{p})^{-1}=R_0
-2\sum_{p_j<p\leq (2{\bar L}(p_j)+1)/D,
(p,D)=1}\bigg[\frac{{\bar L}(p_j)-M(j+1)}{p}
\bigg]\\&&+4\sum_{p_j<p<p'\leq 
(2{\bar L}(p_j)+1)/D,(p,D)=1}\bigg[\frac{
{\bar L}(p_j)-M(j+1)}{pp'}\bigg]\mp\cdots,
\label{tr1}
\end{eqnarray} 
where $[x]$ is the integer part of $x$ as 
usual. Equation~(\ref{tr1}) is equivalent to 
Eq.~(\ref{tr}).~$\diamond$ 

{\bf Definition~5.3.} We split $R=R_M+R_E$ 
into its main and error terms 
\begin{eqnarray}\nonumber
&&R_M\prod_{p|D}(1-\frac{1}{p})^{-1}=R_0
-2\sum_{p_j<p\leq (2{\bar L}
(p_j)+1)/D,(p,D)=1}\frac{{\bar L}(p_j)-
M(j+1)}{p}\\&&+4\sum_{p_j<p<p'\leq 
(2{\bar L}(p_j)+1)/D,(p,D)=1}\frac{
{\bar L}(p_j)-M(j+1)}{pp'}\mp\cdots,\\&&
R_E=2\sum_{p_j<p\leq (2{\bar L}
(p_j)+1)/D,(p,D)=1}\bigg\{\prod_{q|D}(1
-\frac{1}{q})\frac{{\bar L}(p_j)-
M(j+1)}{p}\bigg\}\\&&-4\sum_{p_j<p<p'\leq 
(2{\bar L}(p_j)+1)/D,(p,D)=1}\bigg\{
\prod_{q|D}(1-\frac{1}{q})\frac{{\bar 
L}(p_j)-M(j+1)}{pp'}\bigg\}\mp\cdots,
\label{tr2}
\end{eqnarray}  
using the usual decomposition $[x]=x-\{x\}.$ 

{\bf Theorem~5.4.} {\it The main term $R_M$ 
satisfies}
\begin{eqnarray}\nonumber
&&R_M\prod_{p|D}(1-\frac{1}{p})^{-1}=
{\bar L}(p_j)\prod_{5\leq p\leq 
(2{\bar L}(p_j)+1)/D,(p,D)=1}\left(1-\frac{2}
{p}\right)\\&&+M(j+1)[1-\prod_{p_j<p\leq 
(2{\bar L}(p_j)+1)/D,(p,D)=1}(1-\frac{2}{p})]. 
\label{tr3}
\end{eqnarray}

{\it Proof.} Expanding the product  
\begin{eqnarray}
\prod_{5\leq p\leq p_j,(p,D)=1=(p_j,D)}
(1-\frac{2}{p})
\end{eqnarray}
and combining corresponding sums in 
Eq.~(\ref{tr2})
\begin{eqnarray}\nonumber
&-&\sum_{5\leq p\leq p_j,(p,D)=1=(p_j,D)}\frac{1}
{p}-\sum_{p_j<p\leq (2{\bar L}(p_j)+1)/D,
(p,D)=1}\frac{1}{p}\\&=&-\sum_{5\leq p\leq (2{
\bar L}(p_j)+1)/D,(p,D)=1}\frac{1}{p},\ldots 
\end{eqnarray} 
shifts the upper limit of the primes in 
the product $\prod_p(1-2/p)$ from $p_j$ to 
$(2{\bar L}(p_j)+1)/D$ so that we obtain 
Eq.~(\ref{tr3}). The considerable 
cancellations involved collapse $R_0$ to 
the expected magnitude in $R_M.~\diamond$ 

{\bf Theorem~5.5.} {\it The main term $R_M$ 
obeys the asymptotic law}    
\begin{eqnarray}
R_M\sim \frac{\prod_{p|D}(1-\frac{1}{p})6c_2 
e^{-2\gamma}2{\bar L}(p_j)}{\prod_{p|D}\left(1
-\frac{2}{p}\right)\log^2(2{\bar L}(p_j)+1)/D},
\end{eqnarray}
{\it for $p_j\to\infty,$ where} ${\bar L}(p_j)=
\prod_{5\leq p\leq p_j,(p,D)=1}p.$

{\it Proof.} This follows from Theorem~5.4 as in 
the proof of Theor.~5.8 in 
Ref.~\cite{adhjw}.~$\diamond$ 

\section{Summary and Discussion} 

The twin prime sieves constructed here are  
genuine asymptotic pair sieves that work only 
for prime twins at odd half-distance $D\geq 3.$ 

Accurate counting of non-rank sets require 
the infinite, but sparse set of odd `primorials' 
$\{{\bar L}(p_j)=\prod_{3<p\leq p_j}p,(p,D)=1\}.$ 
The twin primes are not directly sieved, rather 
twin-D-I ranks $m$ are with $2m\pm D$ both prime. 
All other numbers are non-ranks. Primes serve to 
organize and classify (nontrivial) non-ranks in 
arithmetic progressions with equal distances 
(periods) that are primes or products of them. 
     
The coefficients of the asymptotic twin-D-I 
prime distributions depend on $D.$ They are  
$\approx 5$ for sexy primes with $D=3,$  
$\approx 3.3$ for $D=5,$ and $\approx 2.49$ 
for (sufficiently) large $D=$prime for which 
$\prod_{p|D}(1-\frac{1}{p})(1-\frac{2}{p})^{-1}
=1+\varepsilon,$ for some $\varepsilon>0.$ This 
is about a factor $3$ larger than for ordinary 
twins and, remarkably, reflects the different 
abundances allowed by their progressions 
$2m\pm D$ and $2(3m)\pm 1$ in class I. If the 
distance $2D=2\prod_{3\leq p\leq p_e}p,$ then 
the coefficient grows as $\log p_e.$ 

Thus, pair sieves as a resolution of the 
parity problem for prime twins in class I 
allow replacing the need for a lower bound 
on the number of twin-D-I ranks $R$ (or 
$\pi_2/2$) by an upper bound for the remainder 
$R_E$ (that must be lower than $R_M$).


\end{document}